\documentclass{amsart}

\usepackage{graphicx}
\usepackage{amssymb}
\usepackage{latexsym}
\usepackage{color}

\usepackage{amsmath,amsthm}

\usepackage[all]{xy}
\usepackage{mathrsfs}

\usepackage{amsfonts}
\usepackage{amsmath}

\usepackage{mathtools}
\usepackage{dsfont}
\usepackage{url}
\usepackage{enumerate}

\usepackage[hyperindex=true, colorlinks=false]{hyperref}

\newtheorem{theorem}{Theorem}
\newtheorem{conjecture}[theorem]{Conjecture}

\newtheorem{lemma}{Lemma}
\newtheorem{proposition}[lemma]{Proposition}
\newtheorem{corollary}[lemma]{Corollary}

\theoremstyle{definition}

\newtheorem{definition}[lemma]{Definition}
\newtheorem{definition-lemma}[lemma]{Definition-Lemma}
\newtheorem{definition-theorem}[lemma]{Definition-Theorem}
\newtheorem{example}[lemma]{Example}
\newtheorem{remark}[lemma]{Remark}

\newtheorem*{ack}{Acknowledgements}

\newcommand{\HH}{\mathrm{H}}
\newcommand{\zid}{\mathrm{id}}
\newcommand{\BT}{\mathrm{BT}}

\newcommand{\N}{\mathbb{N}}

\newcommand{\R}{\mathbb{R}}
\newcommand{\Z}{\mathbb{Z}}

\newcommand{\D}{\mathbb{D}}
\newcommand{\Sph}{\mathbb{S}}

\newcommand{\fF}{\mathcal{F}}

\newcommand{\fO}{\mathcal{O}}

\newcommand{\pa}{\partial}

\newcommand{\bt}{\operatorname{bt}}

\newcommand{\DSym}{\operatorname{DSym}}

\newcommand{\Fix}{\operatorname{Fix}}

\newcommand{\Homeo}{\operatorname{Homeo}}

\newcommand{\inte}{\operatorname{int}}

\newcommand{\LK}{\operatorname{L}}
\newcommand{\lk}{\operatorname{lk}}
\newcommand{\M}{\operatorname{M}}
\newcommand{\MCG}{\operatorname{MCG}}
\newcommand{\p}{\operatorname{P}}

\newcommand{\Sym}{\operatorname{Sym}}
\newcommand{\Tr}{\operatorname{Tr}}

\newcommand{\GL}{\operatorname{GL}}

\newcommand{\ie}{{\it i.e.}\ }
\newcommand{\st}{{\it s.t.}\ }

\begin{document}


\title{Braids and linked periodic orbits of disc homeomorphisms}

\author{Xiang LIU}
\address{Institute of Mathematics, Academy of Mathematics and System Science, Chinese Academy of Sciences, Beijing 100190, China}%
\email{nz\_liu1989@163.com}

\thanks{Supported by NSF of China (No. 12071309)}%



\begin{abstract}
This paper concerns on linked periodic orbits of orientation-preserving homeomorphisms of the $2$-disc in the sense of Gambaudo. We interpret the linking of periodic orbits by using their induced braids. Then based on the forcing relation of braids, we present a method for finding periodic orbits which are linked with a given one. In particular, for $3$-periodic orbits of pseudo-Anosov braid types, we provide new examples of linked orbits of periods at most $4$.
\end{abstract}

\maketitle



{\it To the memory of my grandfather Zhong-An Liu, Oct 26, 1927-- Sep 23, 2022.}

\section{Introduction}

The geometric structure of periodic orbits is a main subject in both continuous and discrete dynamical systems. For example, it is well studied for (oriented) links composed of closed orbits of $3$-dimensional flows \cite{FS}. A naive fact is that two closed orbits having non-zero linking number are linked.

For a self-map $f\colon X\to X$, a subset $A\subset X$ is {\it invariant} if $f(A)\subset A$. The simplest invariant sets are orbits. The {\it orbit} of a point $x\in X$ is the set
\begin{equation*}
\fO_f (x)= \{f^k (x)\big|k\in \N\}.
\end{equation*}
$x$ is a {\it periodic point} of {\it (least) period} $n$ if its orbit has exactly $n$ elements. In particular, a {\it fixed point} is a $1$-periodic point. The sets of $n$-periodic points and fixed points of $f$ are denoted by $\p_n(f)$ and $\Fix (f)$, respectively.

This paper studies discrete systems in dimension $2$, especially the {\it linking} of periodic orbits of orientation-preserving homeomorphisms of the disc, in the sense introduced by Jean-Marc Gambaudo in 1990 \cite{Gam}. This notion in fact arose from the study of the structure of periodic orbits of orientation-preserving disc diffeomorphisms with topological entropy zero \cite{GamvSS}, where the unlinked was called {\it disjoint}.

Let $M$ be a surface and $h\colon M\to M$ be a homeomorphism.
\begin{definition}
Two invariant sets $A_1, A_2$ of $h$ are called {\it unlinked} if there are disjoint closed discs $D_1, D_2$ in $M$, such that $A_i\subset \inte(D_i)$, and $h(\pa D_i)$ is isotopic to $\pa D_i$ in $M-A_1\cup A_2$, for $i=1,2$.

Otherwise, they are {\it linked}.
\end{definition}
For example, any two fixed points are always unlinked; For a rotation of the disc, any periodic orbit is linked with the center. It is clear that the linking property is invariant under isotopy.

The known results of linked orbits are mainly about a periodic orbit linked with a fixed point, which we called the {\it simple} case. Gambaudo showed that any periodic orbit of an orientation-preserving $C^1$-embedding of $\D$ to itself is linked with a fixed point \cite{Gam}, while Kolev proved the same result for orientation-preserving homeomorphisms of $\R^2$ \cite{K1}. Then Bonino showed that any periodic orbit of period at least $3$ of an orientation-reversing homeomorphism of $\Sph^2$ is linked with an orbit of period $2$ \cite{Bon}. Boro\'nski recently generalized Gambaudo-Kolev Theorem to the case of bounded orbits linked with fixed points of orientation-preserving homeomorphisms of $\R^2$ \cite{Bor}.


The simple case is related to the {\it linking number problem} asked by John Franks, which asserts that for any orientation-preserving homeomorphism $f$ of the plane, any $n$-periodic orbit always has non-zero linking number around a certain fixed point. Bonatti and Kolev confirmed it provided by $f-\zid$ contracting \cite{BK}, while Guaschi proved it for $n=3,4$ \cite{Gu2}. Finally, Le Calvez solved the problem in general \cite{LC}. We note that if the linking number of $p\in\p_n (f)$ around $q\in\Fix (f)$ is non-zero, then the orbit $\fO_f (p)$ is linked with $q$ in the sense of Gambaudo.

Here are two natural problems in this subject. Firstly, are there more {\it non-simple} examples of linked periodic orbits, \ie {\it none} of the two orbits are fixed points? And more fundamentally, what is the mechanism of linking of two orbits?


This paper is motivated by the first one and is inspired by \cite{Gu2}. Focusing on the simplest case, \ie for orientation-preserving homeomorphisms of $\D$, we reduce the linking of orbits to the linking of closed braids they induced in mapping tori, and solve the first problem based on the algorithm of braid forcing proposed in \cite{JZh}.

The phenomenon that a given periodic orbit implies the existence of other orbits is usually called {\it forcing} by the given one. The closely related notions are the {\it braid type} and the {\it forcing relation} of braid types introduced by Boyland and Matsuoka \cite{Boy2, M1}. For an orientation-preserving homeomorphism of $\D$, a periodic orbit induces a braid, whose conjugacy class in the braid group is called the braid type. The braid (type) forces its extensions, some of whom correspond to orbits linked with the given periodic orbit.

Our main result is the following.
\begin{theorem} \label{alg}
Given a set $P$ of $n$ points in $\inte(\D)$, a cyclic braid $\beta\in B_n$, and an integer $m>0$, there is an algorithm to decide a homeomorphism $f\colon\D\to \D$ satisfying that $f|_{\pa\D} =\zid$, $P$ as an $n$-periodic orbit of $f$ induces the braid $\beta$, and $f$ has an $m$-periodic orbit linked with $P$.
\end{theorem}

Therefore, for orientation-preserving disc homeomorphisms, the problem of finding periodic orbits linked with a given one is solved, at least in principle. In particular, to get new non-simple examples of linked periodic orbits, we consider the "simplest" pseudo-Anosov braid type.
\begin{proposition} \label{lpa}
Let $f\colon\D\to \D$ be an orientation-preserving homeomorphism. If $f$ has a $3$-periodic orbit $P$ whose braid type is $[\sigma_1\sigma_2^{-1}]$, then there are $m$-periodic orbits of $f$ linked with $P$ for $m\leq 4$.
\end{proposition}


This paper is organized as follows. In Section 2, we discuss the relationships between periodic orbits, braids, and mapping classes. Then we give a method for finding linked periodic orbits of disc homeomorphisms based on the forcing relation of braids in Section 3. Next in Section 4 we present examples of linked periodic orbits for a pseudo-Anosov $3$-braid. Section 5 devotes further topics.


\section{Periodic orbits, braids, and mapping classes}   \label{secBCHD}

In this section, we review the relationship between periodic orbits and braids to reduce the linked orbits problem, as well as the relationship between braids and mapping classes. We also recall the definition and computation of linking number of two closed braids, which is a part of Theorem~\ref{alg}.


Let $\D$ be the unit closed disc in the plane, and denote the group of orientation-preserving homeomorphisms of $\D$ by $\Homeo_+ (\D)$. By the Alexander trick, any map $f\in\Homeo_+ (\D)$ is isotopic to the identity of $\D$. We fix an isotopy $f_t\colon\zid_\D\sim f$.

For $n\geq 2$ and an $n$-periodic orbit $P$ of $f$, the subset
\begin{equation*}
\beta(P,f_t)= \{(f_t(P),t)\big| t\in I\}
\end{equation*}
of $\D\times I$ is a (geometric) {\it braid} of $n$ strands. We usually omit the isotopy $f_t$ in braid notations. The isotopy classes of $n$-braids form the Artin braid group $B_n$, which has the standard presentation
\begin{equation*}
\big< \sigma_1,\cdots,\sigma_{n-1}\big| \sigma_i\sigma_j=\sigma_j\sigma_i\;\mbox{for}\;\big|i-j\big|\geq 2, \sigma_i\sigma_{i+1}\sigma_i= \sigma_{i+1}\sigma_i\sigma_{i+1}\;\mbox{for}\;1\leq i\leq n-2\big>.
\end{equation*}
There is a canonical group homomorphism $\pi\colon B_n\to S_n$, which maps $\sigma_i$ to the transition $(i,i+1)$ for $i=1,\cdots,n-1$. The braid induced by a periodic orbit is called the {\it braid} of the orbit, and is {\it cyclic} as $\pi$ maps it to an $n$-cycle in $S_n$. More on braid groups see \cite{Bir}.

\begin{center}
\setlength{\unitlength}{1.0mm}
\begin{picture}(20,28)(-20,-8)
\put(0,8){\line(-1,1){16}}
\put(-16,8){\line(1,1){7}}
\put(0,24){\line(-1,-1){7}}
\put(0,-8){\line(1,1){16}}
\put(16,-8){\line(-1,1){7}}
\put(0,8){\line(1,-1){7}}
\put(-16,8){\line(0,-1){16}}
\put(16,8){\line(0,1){16}}

\put(-28,24){\makebox(0,0)[lc]{$t=0$}}
\put(-28,-8){\makebox(0,0)[lc]{$t=1$}}
\end{picture}

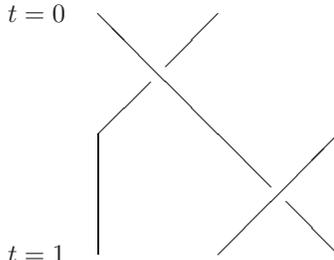
Fig 1.\ The braid $\sigma_1\sigma_2^{-1}\in B_3$
\end{center}

In general, an invariant set $X$ of $N$ points induces an $N$-braid $\beta(X)$, called the {\it braid} of $X$, and periodic orbits in $X$ are bijective correspondence to cycles in $\pi(\beta(X))$. If $P\subset X$ is a periodic orbit, then its braid $\beta(P)$ is a {\it sub-braid} of $\beta(X)$, and we also call $\beta(X)$ an {\it extension} of $\beta(P)$.


The mapping torus is useful in fixed point theory and in our subject. For $f\in\Homeo_+ (\D)$, its {\it mapping torus} is defined to be the quotient space
\begin{equation*}
\M_f= \D\times I\big/ (x,1)\sim (f(x),0).
\end{equation*}
Since $f$ is isotopic to $\zid_\D$, $\M_f$ is a solid torus, and moreover, it is a disc bundle over $\Sph^1$. We denote the corresponding point of $(x,t)$ in $\M_f$ by $[x,t]$. If $f$ is of class $C^1$, there is a canonical flow $\phi_t$ defined on $\M_f$ whose vector field is descended from $\pa_t$ on $\D\times I$. Then the Poincar\'{e} map of $\phi_t$ at the slice $[\D\times 0]$ can be identified with $f\colon\D\to\D$, and $(\M_f,\phi_t)$ is the suspension flow of $f$. We also use the similar construction for orientation-preserving homeomorphisms of punctured discs.

We always suppose that $\M_f$ is embedded in $\Sph^3$ trivially. For an $n$-periodic orbit $P$, the {\it closure} $\hat{\beta}(P)$ of the braid $\beta(P)$ is obtained by identifying the initial points and end points of each of the strands. It is a simple closed curve contained in $\M_f$ and hence embedded in $\Sph^3$. We orient it in the direction of increasing $t$, thus $\hat{\beta}(P)$ is an oriented knot.

Here is the key observation.
\begin{proposition} \label{iff}
Two periodic orbits $P$ and $Q$ of $f\in\Homeo_+(\D)$ are linked if and only if the corresponding closed braids $\hat{\beta}(P)$ and $\hat{\beta}(Q)$ are linked in $\M_f(\subset\Sph^3)$.
\end{proposition}
The only if part is similar to the Corollary in Gambaudo's original paper~\cite{Gam}, and the if part comes from Guaschi's paper~\cite{Gu2}.


We determine whether two such closed braids are linked by using the linking number. Two closed braids are linked if they have non-zero linking number, but the converse is not true. Firstly, we can read the linking number geometrically from the diagram of the braids. In a diagram of two cyclic braids $\beta,\alpha$, for each crossing $C$ where the arc of $\beta$ crosses under the arc of $\alpha$, we assign the {\it sign} $\varepsilon(C)$ of $C$ as in Fig 2.
\begin{center}
\setlength{\unitlength}{1.0mm}
\begin{picture}(20,28)(-10,-14)
\put(-8,-8){\line(1,1){16}}
\put(8,-8){\line(-1,1){7}}
\put(-8,8){\line(1,-1){7}}

\put(-8,-8){\line(1,0){1}}
\put(-8,-8){\line(0,1){1}}
\put(8,-8){\line(-1,0){1}}
\put(8,-8){\line(0,1){1}}

\put(8,0){\makebox(0,0)[lc]{$-1$}}
\put(-12,8){\makebox(0,0)[lc]{$\beta$}}
\put(10,8){\makebox(0,0)[lc]{$\alpha$}}
\end{picture}
\ \ \ \ \ \ \ \ \ \ \ \ \ \ \ \
\begin{picture}(20,28)(-10,-14)
\put(-8,-8){\line(1,1){7}}
\put(8,8){\line(-1,-1){7}}
\put(-8,8){\line(1,-1){16}}

\put(-8,-8){\line(1,0){1}}
\put(-8,-8){\line(0,1){1}}
\put(8,-8){\line(-1,0){1}}
\put(8,-8){\line(0,1){1}}

\put(8,0){\makebox(0,0)[lc]{$+1$}}
\put(-12,8){\makebox(0,0)[lc]{$\alpha$}}
\put(10,8){\makebox(0,0)[lc]{$\beta$}}
\end{picture}

Fig 2.\ The sign of a crossing
\end{center}
Then the linking number of the two closed braids satisfies
\begin{equation*}
\lk(\hat{\beta},\hat{\alpha})=\sum_C \varepsilon(C).
\end{equation*}

Next we recall the homological definition of linking number of two periodic orbits and a computation method due to John Guaschi \cite{Gu2} by using representations of braid groups.

For $f\in\Homeo_+(\D)$ with an $n$-periodic orbit $P$, fix an isotopy $f_t\colon\zid_\D \sim f$. Let $\{\gamma_1,\gamma_2\}$ be the standard basis of
$\HH_1 (\M_{f|_{\D-P}})\cong\Z^2$,
where $\gamma_1$ corresponds to an oriented meridian of the closed braids $\hat{\beta}(P)$, and $\gamma_2$ corresponds to a longitude on
$\pa (\M_{f|_{\D-P}})$ oriented by increasing $t$. Then if $Q$ is an $m$-periodic orbit of $f$, the loop $\{[Q\times t]\big|t\in I\}$ is a simple closed curve in $\M_{f|_{\D-P}}$, which gives a homology class
\begin{equation*}
l\gamma_1+m\gamma_2\in\HH_1 (\M_{f|_{\D-P}}).
\end{equation*}
We call $l$ the {\it linking number $\LK(Q,P)$ of $Q$ about $P$ (relative to $f_t$)}. We note that for $m=1$, $\LK(q,P)=\LK_{f,p}(q)$ the linking number of $P$ around a fixed point $q$ with $p\in P$, the linking number $\LK$ is symmetric as
\begin{equation*}
\LK(Q,P)=\lk(\hat{\beta}(Q),\hat{\beta}(P))=\lk(\hat{\beta}(P),\hat{\beta}(Q))=\LK(P,Q),
\end{equation*}
and $P,Q$ are linked if $\LK(Q,P)\neq 0$. See more details in \cite{Gu2}.


Artin in 1925 showed that $B_n$ acts faithfully on the free group $F_n=\big<x_1,\cdots,x_n\big>$ from right by
\begin{equation*}
x_i\cdot\sigma_i=x_ix_{i+1}x_i^{-1}, \quad x_{i+1}\cdot\sigma_i=x_i, \quad x_j\cdot\sigma_i=x_j
\end{equation*}
for $j\neq i,i+1$.
For a braid $\beta\in B_n$, suppose that its closure $\hat{\beta}$ is a link with $\mu\geq 1$ components $L_1,\cdots,L_{\mu}$. Then the link group $G=\pi_1(\Sph^3-\hat{\beta})$ has the presentation
\begin{equation*}
\big<x_1,\cdots,x_n\big|R_1,\cdots,R_n\big>,
\end{equation*}
where $R_i=A_i x_{\rho(i)}A_i^{-1}x_i^{-1}$ for some $\rho\in S_n$ consisting of disjoint cycles $\rho_j$ corresponding to $L_j$ and $A_i\in F_n$, $i=1,\cdots,n$.
There is a canonical homomorphism
\begin{equation*}
\Psi\colon\Z F_n\to\Z G
\end{equation*}
of group rings induced by the quotient homomorphism $F_n\to G$.

We consider a free Abelian group $H=\big<t_1\big>\times\cdots\times\big<t_\mu\big>$, where $t_j$ corresponds to the component $L_j$.
Then we define a homomorphism
\begin{equation*}
\Phi\colon\Z G\to\Z H, \quad x_i\mapsto t_j,
\end{equation*}
whenever $i$ is in the cycle $\rho_j$. We note that the group ring $\Z H=\Z [t_1^{\pm 1},\cdots,t_\mu^{\pm 1}]$.
Now let $g_i=x_1\cdots x_i$ for $i=1,\cdots,n$. Then $\{g_1,\cdots,g_n\}$ freely generates $F_n$ as $x_i=g_{i-1}^{-1}g_i$, and $g_n\cdot\beta=g_n$ for any $\beta\in B_n$.

The {\it Magnus representation} $R\colon B_n\to \GL (n,\Z H)$ of $B_n$ is defined by
\begin{equation*}
R(\beta)=\Big(\Phi\circ\Psi\Big(\frac{\pa (g_i\cdot\beta)}{\pa g_j}\Big)\Big)_{1\leq i,j\leq n},
\end{equation*}
where $\pa\big/\pa g_j$ are the {\it Fox's free derivatives} given by
\begin{equation*}
\frac{\pa g_i}{\pa g_j}=\delta_{ij},\quad \frac{\pa (vw)}{\pa g_j}=\frac{\pa v}{\pa g_j}+v\frac{\pa w}{\pa g_j},
\end{equation*}
for $1\leq i,j\leq n$ and words $v,w$ in $\big<g_1,\cdots,g_n\big>$.
Then the {\it link representation}
\begin{equation*}
r\colon B_n\to \GL (n-1,\Z H)
\end{equation*}
is given by
\begin{equation*}
R(\beta)=
\begin{bmatrix}
r(\beta) & * \\
0\cdots 0 & 1
\end{bmatrix}.
\end{equation*}

Now we state Guaschi's formula for computing the linking number of two closed braids. Here $\Z H=\Z [t_1^{\pm 1},t_2^{\pm 1}]$.
\begin{theorem}[\cite{Gu2}, for $\mu =2$] \label{lk}
Let $\beta\in B_n$ be a sub-braid of $\beta\cup\alpha\in B_{n+m}$, whose closure is a link with two components $\hat{\beta}, \hat{\alpha}$. Then
\begin{equation*}
\det(r(\beta\cup\alpha)-I)\big|_{t_2=1} =(-1)^{m-1}(t_1^l-1)\det(r(\beta)-I),
\end{equation*}
where $l=\lk(\hat{\beta}, \hat{\alpha})$.
\end{theorem}

\begin{remark} \label{gln}
These representations of the braid group are also useful to compute the generalized Lefschetz number as
\begin{equation*}
L_H(\beta)= 1-\Tr(R(\beta))= -\Tr(r(\beta))\in\Z H.
\end{equation*}
\end{remark}


\begin{example} \label{lkn}
Let $\beta=\sigma_1\sigma_2^{-1}\in B_3$. Then $\det(r(\beta)-I)= 1+t_1+t_1^{-1}$.

For $m=2$, $\beta$ is a sub-braid of $\beta\cup\alpha_1= \sigma_1\sigma_2^{-1}\sigma_3^{-2}\sigma_4^{-1}\in B_5$. Then we have $\lk(\hat{\beta},\hat{\alpha}_1)=-1$, since
\begin{equation*}
\det(r(\beta\cup\alpha_1)-I)\big|_{t_2=1} = -(t_1^{-1}-1)(1+t_1+t_1^{-1}).
\end{equation*}
Similarly, $\beta$ is also a sub-braid of $\beta\cup\alpha_2= \sigma_1\sigma_2^{-1}\sigma_4\sigma_3\sigma_2\sigma_1^2\sigma_2^{-1}\sigma_3^{-1}\in B_5$, and $\lk(\hat{\beta},\hat{\alpha}_2)=1$.

For $m=3$, $\beta$ is a sub-braid of
\begin{equation*}
\beta\cup\gamma_1= \sigma_1\sigma_2^{-1}\sigma_3\sigma_4\sigma_5\sigma_2\sigma_3\sigma_4\sigma_1\sigma_2\sigma_3\sigma_1\sigma_2\sigma_3\sigma_4^{-1}\sigma_5^{-1}\sigma_1\sigma_2\sigma_3^{-1}\sigma_4^{-1}\sigma_1\sigma_2^{-1}\sigma_3^{-3}\sigma_4^{-1}
\in B_6,
\end{equation*}
and $\lk(\hat{\beta},\hat{\gamma}_1)=2$.

For $m=4$, $\beta$ is a sub-braid of
\begin{equation*}
\begin{array}{rl}
\beta\cup\delta_1 = & \sigma_1\sigma_2^{-1}\sigma_4\sigma_5\sigma_6\sigma_3\sigma_4\sigma_5\sigma_2\sigma_3\sigma_4\sigma_1\sigma_2\sigma_3\sigma_1\sigma_2\sigma_1\sigma_3\sigma_2\sigma_1 \\
  & \cdot\sigma_4^{-1}\sigma_3^{-1}\sigma_2^{-1}\sigma_5^{-1}\sigma_4^{-1}\sigma_3^{-1}\sigma_6^{-1}\sigma_5^{-1}\sigma_4^{-1}\sigma_3^{-1}\sigma_4^{-1}\sigma_3^{-1}\sigma_4^{-1}\sigma_3^{-1}\sigma_4^{-1}\in B_7,
\end{array}
\end{equation*}
and $\lk(\hat{\beta},\hat{\delta}_1)=1$.
\qed
\end{example}


At the end of this section, we recall the Nielsen-Thurston classification of mapping classes and how braids connect with them. We begin with pseudo-Anosov maps, which are incisively introduced by William Thurston in 1970s.

A homeomorphism $\phi\colon M\to M$ of a surface $M$ is called {\it pseudo-Anosov} if there is a number $\lambda>1$ and a pair of transverse measured foliations $(\fF^u,\mu^u)$ and $(\fF^s,\mu^s)$ on $M$, \st
\begin{equation*}
\phi(\fF^u,\mu^u)=(\fF^u,\lambda\mu^u),\quad \phi(\fF^s,\mu^s)=(\fF^s,\frac{1}{\lambda}\mu^s).
\end{equation*}
We only emphasize that pseudo-Anosov maps have useful dynamical properties to be the candidate for presenting new examples of linked periodic orbits.
\begin{proposition} \label{dpa}
Let $\phi\colon M\to M$ be pseudo-Anosov. Then

(1) the set of periods of $\phi$ is cofinite in $\N$, and

(2) every periodic point of $\phi$ is unremovable.
\end{proposition}
The first property due to Gambaudo and Llibre \cite{GL} guarantees that given a periodic orbit, there are plenty of other periodic orbits. The second one, noted by Asimov and Franks, means that periodic orbits are persistent under isotopy \cite{AF}.

Now we state the Nielsen-Thurston classification theorem \cite{T}, which is fundamental in geometric topology.
\begin{theorem} \label{NT}
Let $f\colon M\to M$ be a homeomorphism of a compact connected surface $M$. Then there is a homeomorphism $\phi$ isotopic to $f$, called its Thurston form, satisfying that
\begin{itemize}
  \item either $\phi$ is of finite type, \ie $\phi^m=\zid$ for some $m\geq 1$; or
  \item $\phi$ is pseudo-Anosov; or
  \item $\phi$ is reducible, \ie there is a system of disjoint simple closed curves $\Gamma=\{\gamma_1,\cdots,\gamma_k\}$ in $\inte(M)$, called reducing curves, \st $\Gamma$ is $\phi$-invariant.
\end{itemize}
Moreover, pseudo-Anosov maps are neither of finite type nor reducible.
\end{theorem}
The isotopy class of $f$ is called its {\it mapping class}, thus we have the Nielsen-Thurston classification for mapping classes. Mapping classes form a group with multiplication induced by the composition of maps, which is called the {\it mapping class group} of $M$ and denoted by $\MCG(M)$. If $M$ is oriented, then $\MCG(M)= \pi_0 (\Homeo_+(M))$.


Braids and mapping classes are naturally related as follows. Given punctures $X=\{p_1,\cdots,p_n\}$ in $\inte(\D)$ and a braid $\beta\in B_n$, sliding $\D$ down along $\beta$ gives an orientation-preserving homeomorphism $h\colon\D\to\D$ satisfying $h(X)=X$. We denote the subgroup of such homeomorphisms by $\Homeo_+(\D,X)$. If, in addition, we require that $\pa\D$ is point-wise fixed, the corresponding mapping class group is denoted by $\MCG(\D,X;\pa\D)$. Furthermore, we have the following
\begin{proposition}[\cite{Bir}] \label{iso}
There is a canonical isomorphism
\begin{equation*}
\lambda\colon B_n\to\MCG(\D,X;\pa\D)
\end{equation*}
mapping $\sigma_i$ to the twist exchanging $p_i,p_{i+1}$.
\end{proposition}
The above sliding down the disc along a braid exactly realizes the isomorphism $\lambda$, after perturbing the homeomorphism to fix the boundary point-wise.

Under $\lambda$, the conjugacy class of a braid in $B_n$ corresponds to the conjugacy class of its image in $\MCG(\D,X;\pa\D)$, and the former is exactly the braid type we consider in the next section. On the other hand, by the Nielsen-Thurston classification, braids and their conjugacy classes naturally fall into three types induced by $\lambda$: finite type, pseudo-Anosov, and reducible.


\section{Forcing of braids and the algorithm of linked periodic orbits}

In this section, we recall the notion of braid types, which focuses on the existence of periodic orbits implied by a given orbit. Then relying on the algorithm of forced extensions due to Jiang and Zheng \cite{JZh}, we present our method for finding linked periodic orbits of disc homeomorphisms, which is the proof of Theorem~\ref{alg}.


The notion of braid types was introduced by Philip Boyland and Takashi Matsuoka in 1980s for characterizing the forcing relation in dimension $2$, in analog to the role of periods in the Sharkovsky ordering for interval self-maps \cite{M1, Boy2}.
\begin{definition}
For $f\in \Homeo_+(\D)$ with an $n$-periodic orbit $P$ in $\inte(\D)$, the {\it braid type} $\bt(P,f)$ of $P$ is defined to be the conjugacy class of the braid $\beta(P,f_t)$ in $B_n$ for an isotopy $f_t\colon\zid_\D\sim f$.
\end{definition}
It is easy to see that $\bt(P,f)$ is independent of the choice of the isotopy $f_t$. Thus an $n$-braid type is in fact a conjugacy class in the braid group $B_n$. It is equivalent to define $n$-braid types to be the conjugacy classes in $\MCG(\D,P;\pa\D)$ for any subset $P\subset\inte(\D)$ of $n$ points by the isomorphism $\lambda$, or to be the free loop classes in the unordered configuration space of $n$ points in $\D$ \cite{M2}. Let
\begin{equation*}
\BT_n=\{\mbox{braid types of}\; n\mbox{-periodic orbits for all homeomorphisms of}\; \D\},
\end{equation*}
and let $\BT=\bigcup_{n\geq 1}\BT_n$ be the set of all braid types.

\begin{definition}
We say $f\in\Homeo_+(\D)$ {\it exhibits} a braid type $[\beta]\in\BT_n$ if $f$ has an $n$-periodic orbit $P$ whose braid type is $[\beta]$. Then we say a braid type $[\beta_1]\in\BT$ {\it forces} $[\beta_2]\in\BT$, denoted by $[\beta_1]\geq[\beta_2]$, if once $f\in\Homeo_+(\D)$ exhibits $[\beta_1]$ then $f$ also exhibits $[\beta_2]$.
\end{definition}
We also say the forcing of braids, \ie a braid $\beta_1$ {\it forces} another braid $\beta_2$ if $[\beta_1]\geq[\beta_2]$. It turns out that the forcing relation $\geq$ is a partial order on $\BT$ \cite{Boy1}.


To reveal the extra information about how a forced braid winds around the given one, BoJu Jiang and Hao Zheng considered the forcing problem on extensions $\beta(P\cup Q)$ of a given braid $\beta(P)$, in stead of just on $\beta(Q)$. We apply their criterion for forced extensions in the problem on linked periodic orbits, as the most significant step in our algorithm.

Here are their main notions and results. Suppose that $\beta'\in B_{n+m}$ is an extension of $\beta\in B_n$, and let $\phi$ be the Thurston form of mapping class $\lambda(\beta')$.
\begin{definition}
The extension $\beta'$ is {\it collapsible} (resp. {\it peripheral}) relative to $\beta$ if there exists a system of reducing curves of $\phi$, \st one of them encloses none of (resp. precisely one of or all of) the punctures corresponding to $\beta$.
\end{definition}
\begin{proposition} [\cite{JZh}]
An extension $\beta'$ of a pseudo-Anosov braid $\beta$ is collapsible or peripheral relative to $\beta$ if and only if $\beta'$ is reducible.
\end{proposition}

To give the conditions of forced extensions, Jiang and Zheng constructed a group representation
\begin{equation*}
\zeta_{n,m}\colon B_n\to\GL(N,\Z B_{n+m})
\end{equation*}
for $N=\tbinom{m+n-2}{m}$, which is the dual of another matrix representation $\xi_{n,m}$ relating to the Lawrence-Krammer-Bigelow representation \cite{Zh}. Then they got the following trace formula by applying Nielsen fixed point theory on the $m$-th symmetric product space of $Y_n$, the disc with $n$ open discs removed, which is the compactification of $\D-P$. For a nontrivial braid $\beta\in B_n$ and $m\geq 1$,
\begin{equation*}
L_{\Gamma_{\beta,m}} (\Sym^m(\bar f_{\beta}))= (-1)^m \Tr_{\Gamma_{\beta,m}} \zeta_{n,m}(\beta) -\mbox{collapsible terms},
\end{equation*}
where $\bar f_{\beta}\colon Y_n\to Y_n$ is the blowing-up (introduced by Bowen \cite{Bow}) of a representative $f_{\beta}$ of $\lambda(\beta)$, $\Sym^m(\bar f_{\beta})$ is the induced map of $\bar f_{\beta}$ on the symmetric product space, and $\Gamma_{\beta,m}$ is the fundamental group of the mapping torus of the deleted symmetric product map $\DSym^m(\bar f_{\beta})$ being identified with a subgroup of $B_{n+m}$. Note that the generalized Lefschetz number $L_{\Gamma_{\beta,m}} (\Sym^m(\bar f_{\beta}))$ takes value in the free Abelian group generated by conjugacy classes in $\Gamma_{\beta,m}$. Then the $(n+m)$-forced extensions of $\beta$ are exactly the {\it non-peripheral} terms in $L_{\Gamma_{\beta,m}} (\Sym^m(\bar f_{\beta}))$. More details see \cite{JZh}, we only state their criterion here.
\begin{theorem} [\cite{JZh}] \label{ext}
Let $\beta'\in B_{n+m}$ be an extension of $\beta\in B_n$. Then $\beta$ forces $\beta'$ if and only if
\begin{itemize}
  \item $\beta'$ is neither collapsible nor peripheral relative to $\beta$, and
  \item the conjugacy class $[\beta']$ has a non-zero coefficient in $\Tr_{B_{n+m}}\zeta_{n,m}(\beta)$.
\end{itemize}
\end{theorem}


Combining the theorems of Guaschi and Jiang-Zheng, we can now present our algorithm of linked periodic orbits for disc homeomorphisms.

{\it Proof of Theorem~\ref{alg}}.
The inputs of the algorithm are the following: a set $P\subset\inte(\D)$ of $n$ points as the given $n$-periodic orbit; a cyclic $n$-braid $\beta$, which will be the braid of $P$; and an integer $m>0$ as the period of orbits linked with $P$.

1. Slide the disc $\D$ down along $\beta$ to obtain a map $f\in\Homeo_+(\D,P)$, which is a representative of $\lambda(f)$ after perturbing as described in Proposition~\ref{iso}.

2. Apply Jiang-Zheng's Theorem~\ref{ext} to determine forced extensions $\beta\cup\alpha\in B_{n+m}$ of $\beta$. Some of them may have sub-braids $\alpha$ that correspond to periodic orbits linked with $P$.

3. Identify sub-braids $\alpha$ satisfying the condition $\lk(\hat{\beta},\hat{\alpha})\neq 0$, by using Guaschi's formula in Theorem~\ref{lk} algebraically or by reading diagrams geometrically.

4. Locate the periodic orbits $Q$ linked with $P$ to be the endpoints of $\alpha$ obtained from the last step.

Then we get the output including an orientation-preserving homeomorphism $f\colon\D\to\D$ with a given $n$-periodic orbit $P$ of braid $\beta$, and several $m$-periodic orbits $Q$ linked with $P$.
\qed

We give some remarks on the algorithm.
\begin{remark}
If $\beta\in B_n$ is a pseudo-Anosov braid, then for each $m>0$, there are only finitely many braid types $[\beta\cup\alpha]$ forced by $[\beta]$ \cite{JG}. Thus Step 2 yields finite candidates.
\end{remark}
\begin{remark}
In general, it is complicated to compute the trace $\Tr_{B_{n+m}}\zeta_{n,m}(\beta)$ by hand. However, for $m=1$, there is a smaller forcing relation, for which computing of forced $(n+1)$-extensions of a given $n$-braid is easier \cite{WZh}.
\end{remark}

We should note that in Step 3, the case of linked closed braids of linking number zero can be identified geometrically but not algebraically. For an example, see Remark~\ref{lk0} in the next section.


\section{Examples of linked periodic orbits of pseudo-Anosov $3$-braid types}

We give new examples of linked periodic orbits in this section as the proof and corollary of Proposition~\ref{lpa}, by concerning pseudo-Anosov $3$-braid types. As mentioned in Proposition~\ref{dpa}, pseudo-Anosov maps are sufficiently complicated in dynamics to provide linked periodic orbits. On the other hand, the forcing relation among pseudo-Anosov $3$-braid types can be described clearly.

Matsuoka pointed out that pseudo-Anosov $3$-braid types are represented in the $L$-$R$ form \cite{M1, Ha}.
\begin{proposition}
If a $3$-braid type is not of finite type, then it has a representative $w\in B_3$ only consisting of letters $L=\sigma_1$ and $R=\sigma_2^{-1}$.

Moreover, it is pseudo-Anosov if and only if $w$ contains at least one $L$ and one $R$.
\end{proposition}
Therefore, $\sigma_1\sigma_2^{-1}=LR$ represents the {\it simplest} pseudo-Anosov $3$-braid type in the sense of word-length. Furthermore, Michael Handel proved that the forcing relation among pseudo-Anosov $3$-braid types is distinct in terms of the $L$-$R$ form \cite{Ha}.
\begin{theorem} \label{frpa}
For two pseudo-Anosov $3$-braid types $[w],[v]$ in the $L$-$R$ form, $[w]\geq[v]$ if and only if $v$ is obtained from $w$ by removing letters and by cyclic permutation.
\end{theorem}
For example, $LRLLRR$ forces $LRLL,LRLR,LRRR,LR$.


In particular, any pseudo-Anosov $3$-braid type forces $[\sigma_1\sigma_2^{-1}]$. So we start new non-simple examples with $[\sigma_1\sigma_2^{-1}]$ as Proposition~\ref{lpa}.

{\it Proof of Proposition \ref{lpa}}. As mentioned in \cite{Bon}, after perturbing $f|_{\D-P}$ by blowing-up, Gambaudo's theorem holds for homeomorphisms. Thus the case of $m=1$ is trivial. More precisely, as Remark~\ref{gln} the generalized Lefschetz number is
\begin{equation*}
L_H (\sigma_1\sigma_2^{-1})=-1+t_1+t_1^{-1},
\end{equation*}
so $f\colon\D\to\D$ has exactly two fixed points linked with $P$ of non-zero linking numbers $1$ and $-1$, respectively.

For $m=2$, $\beta=\sigma_1\sigma_2^{-1}$ forces an extension $\beta\cup\alpha_1=\sigma_1\sigma_2^{-1}\sigma_3^{-2}\sigma_4^{-1}\in B_5$ provided in Example 6.4 in \cite{JZh} and of linking number $\lk(\hat{\beta},\hat{\alpha}_1)=-1$ as in Example~\ref{lkn}.

For $m=3$, $\beta=\sigma_1\sigma_2^{-1}$ forces an extension
\begin{equation*}
\beta\cup\gamma_1= \sigma_1\sigma_2^{-1}\sigma_3\sigma_4\sigma_5\sigma_2\sigma_3\sigma_4\sigma_1\sigma_2\sigma_3\sigma_1\sigma_2\sigma_3\sigma_4^{-1}\sigma_5^{-1}\sigma_1\sigma_2\sigma_3^{-1}\sigma_4^{-1}\sigma_1\sigma_2^{-1}\sigma_3^{-3}\sigma_4^{-1}\in B_6
\end{equation*}
provided in \cite{J3} and of linking number $\lk(\hat{\beta},\hat{\gamma}_1)=2$ as in Example~\ref{lkn}.

For $m=4$, $\beta=\sigma_1\sigma_2^{-1}$ forces an extension
\begin{equation*}
\begin{array}{rl}
\beta\cup\delta_1 = & \sigma_1\sigma_2^{-1}\sigma_4\sigma_5\sigma_6\sigma_3\sigma_4\sigma_5\sigma_2\sigma_3\sigma_4\sigma_1\sigma_2\sigma_3\sigma_1\sigma_2\sigma_1\sigma_3\sigma_2\sigma_1 \\
  & \cdot\sigma_4^{-1}\sigma_3^{-1}\sigma_2^{-1}\sigma_5^{-1}\sigma_4^{-1}\sigma_3^{-1}\sigma_6^{-1}\sigma_5^{-1}\sigma_4^{-1}\sigma_3^{-1}\sigma_4^{-1}\sigma_3^{-1}\sigma_4^{-1}\sigma_3^{-1}\sigma_4^{-1}\in B_7
\end{array}
\end{equation*}
provided in \cite{J3} and of linking number $\lk(\hat{\beta},\hat{\delta}_1)=1$ as in Example~\ref{lkn}.
\qed

\begin{remark} \label{lk0}
For $m=1$, by using the algorithm of Theorem~\ref{alg}, $\beta=\sigma_1\sigma_2^{-1}$ forces three extensions $\beta A_{1,4}$, $\beta A_{1,4}A_{2,4}^{-1}$, and $\beta A_{2,4}^{-1}$ of linking numbers $1,0$ and $-1$, respectively. But the braid diagram tells us that the closure of the fourth strand of $\beta A_{1,4}A_{2,4}^{-1}$ is indeed linked with $\hat\beta$, so this fixed point is also linked with $P$ by Proposition~\ref{iff}. It corresponds to the term $-1$ of the generalized Lefschetz number.
\end{remark}

According to Theorem~\ref{frpa}, an immediate consequence of Proposition~\ref{lpa} is the following.
\begin{corollary} \label{lpac}
If $f\in\Homeo_+(\D)$ has a $3$-periodic orbit whose braid type is pseudo-Anosov, then for each $m\leq 4$, there are two periodic orbits of periods $3$ and $m$, respectively, which are linked.
\end{corollary}


We want to mention that Proposition~\ref{lpa} can be generalized not only for larger $m$, but also for other braid types. For instance, Eiko Kin studied two families of braids $\beta_{m,n}$ and $\sigma_{m,n}$ in $B_{m+n+1}$ for $m,n\geq 1$, where $\beta_{1,1}=\sigma_1\sigma_2^{-1}$, $\beta_{m,n}$ is pseudo-Anosov, and $\sigma_{m,n}$ is pseudo-Anosov if and only if $\big|m-n\big|\geq 2$. She described the forcing relation among braid types of these two families precisely by using the train track \cite{Ki}. Thus in principle, one can apply the algorithm to these families of braids to provide more examples of linked periodic orbits as in Proposition~\ref{lpa} and Corollary~\ref{lpac}, though the computation is harder than that for pseudo-Anosov $3$-braids.

Here is a general conjecture on linked periodic orbits of pseudo-Anosov braid types.
\begin{conjecture}
If $f\in\Homeo_+(\D)$ has a periodic orbit $P$ whose braid type is pseudo-Anosov, then there is a number $N>0$, \st for each $m>N$, $f$ has {\color{blue}an} $m$-periodic orbit linked with $P$.
\end{conjecture}
We should note that it is still unknown whether each pseudo-Anosov $m$-braid type for $m>3$ is forced by a certain $3$-braid type.


\section{Final remarks}
\label{secMldCalc}

We propose an algorithm of linked periodic orbits of orientation-preserving disc homeomorphisms, relying on Jiang-Zheng's theorem on forced extensions and Guaschi's formula on the linking number of braids. The former is obtained by analyzing fixed point classes of the induced map on the symmetric product space of the punctured disc, while the later comes from translating the linking information of the strands to coordinates of the fixed point classes. Both of them are essentially refinements of {\it Nielsen fixed point theory} \cite{J1}.


It is natural to seek linked periodic orbits on other surfaces, such as the annulus, torus, and hyperbolic surfaces. However, the braid type theory of surface braids is almost an unknown world. Besides, the forcing relation on the annulus is related to rotation sets \cite{Boy2, GarM}.

Here are two facts coming from the {\it template theory} of $3$-dimensional flows. Philip Holmes and R. F. Williams in 1985 pointed out that suspensions of Smale horseshoe have closed orbits which are linked torus knots \cite{HolW}. Later, Holmes showed that any two closed orbits of suspensions of the Josephson equation
\begin{equation*}
\dot{\theta}=r,\quad \dot{r}= -\sin\theta+\nu-\delta r+\beta\cos t,\quad (\theta,r)\in \Sph^1\times I.
\end{equation*}
are linked \cite{Hol}. Consequently, we have the following proposition.
\begin{proposition}
(1) Systems containing Smale horseshoe have non-simple linked periodic orbits.

(2) Any two periodic orbits of the Poincar\'e map of the Josephson equation are linked.
\end{proposition}


If one fixes an area form on a compact surface, then after perturbing one can choose a monotone symplectomorphism as the representative in each mapping class \cite{S}. This leads the {\it symplectic Floer homology}, whose chain complex is generated by non-degenerate fixed points. The structure of the symplectic Floer homology of all the three Nielsen-Thurston types is clear \cite{Gau, CC, N}. More generally, we want to know what can the variants of Floer theory do about linked periodic orbits, including the periodic Floer homology \cite{HS} and the braid Floer homology \cite{VGVW}.


\begin{ack}
The author would like to thank Prof. XueZhi Zhao at Capital Normal University and JianLu Zhang at AMSS for discussions, Jun Wang at Hebei Normal University for the help on references, and Prof. JianZhong Pan and Yang Su at AMSS for their encouragement. The author also would like to thank the referee for carefully reading the manuscript and helpful comments.
\end{ack}


\end{document}